\documentclass[12pt]{amsart}
\usepackage{amsmath, amsfonts, amssymb, amsthm,hyperref}
\hypersetup{hypertex=true,
	colorlinks=true,
	linkcolor=blue,
	anchorcolor=blue,
	citecolor=blue}
\usepackage{bm}
\allowdisplaybreaks[4]
\def\({\bg(}
\def\){\bg)}

\def\v{{\bm v}}

\def\pmod #1{\ ({\rm{mod}}\ #1)}

\theoremstyle{plain}
\newtheorem{theorem}{Theorem}[section]
\newtheorem{lemma}{Lemma}

\newtheorem{conjecture}{Conjecture}

\theoremstyle{definition}

\theoremstyle{remark}
\newtheorem{remark}{Remark}

\makeatletter
\@namedef{subjclassname@2020}{%
	\textup{2020} Mathematics Subject Classification}
\makeatother
\vspace{4mm}

\begin{document}
	\medskip
	
	\title[matrices concerning multiplicative subgroups of finite fields]
	{A conjecture of Zhi-Wei Sun on matrices concerning multiplicative subgroups of finite fields}
	\author[J. L and H.-L. Wu]{Jie Li and Hai-Liang Wu*}

	\address {(Jie Li) School of Science, Nanjing University of Posts and Telecommunications, Nanjing 210023, People's Republic of China}
	\email{\tt lijiemath@163.com}
	
	\address {(Hai-Liang Wu) School of Science, Nanjing University of Posts and Telecommunications, Nanjing 210023, People's Republic of China}
	\email{\tt whl.math@smail.nju.edu.cn}
	
	\keywords{Legendre symbols, Finite Fields, Cyclotomic Matrices, Determinants.
		\newline \indent This work was supported by the Natural Science Foundation of China (Grant No. 12101321).	
        \newline \indent *Corresponding author.}
	
	\begin{abstract}
	Motivated by the recent work of Zhi-Wei Sun on determinants involving the Legendre symbol, in this paper, we study some matrices concerning subgroups of finite fields. 
	
	For example, let $q\equiv 3\pmod 4$ be an odd prime power and let $\phi$ be the unique quadratic multiplicative character of the finite field $\mathbb{F}_q$. If set $\{s_1,\cdots,s_{(q-1)/2}\}=\{x^2:\ x\in\mathbb{F}_q\setminus\{0\}\}$, then we prove that 
	$$\det\left[t+\phi(s_i+s_j)+\phi(s_i-s_j)\right]_{1\le i,j\le (q-1)/2}=\left(\frac{q-1}{2}t-1\right)q^{\frac{q-3}{4}}.$$
	This confirms a conjecture of Zhi-Wei Sun. 
	\end{abstract}
	\maketitle
	
	\section{Introduction}
	\setcounter{lemma}{0}
	\setcounter{theorem}{0}
	\setcounter{equation}{0}
	\setcounter{conjecture}{0}
	\setcounter{remark}{0}
	\setcounter{corollary}{0}

    Let $p$ be an odd prime. The research of determinants involving the Legendre symbol $(\frac{\cdot}{p})$ can be traced back to the works of Lehmer \cite{Lehmer}, Carlitz \cite{Carlitz} and Chapman \cite{Chapman}. For example, Carlitz \cite[Theorem 4]{Carlitz} studied the determinant
    \begin{equation*}
    	\det C(t):=\det\left[t+\left(\frac{i-j}{p}\right)\right]_{1\le i,j\le p-1}.
    \end{equation*}
    Carlitz showed that 
    $$\det C(t)=(-1)^{\frac{p-1}{2}}p^{\frac{p-3}{2}}\left((p-1)t+(-1)^{\frac{p-1}{2}}\right).$$
    
    Along this line, Chapman \cite{Chapman} further investigated some variants of $\det C(t)$. For instance, Chapman considered 
    $$\det C_1(t):=\det \left[t+\left(\frac{i+j-1}{p}\right)\right]_{1\le i,j\le (p-1)/2}.$$
    If we let $\varepsilon_p>1$ and $h_p$ be the fundamental unit and the class number of $\mathbb{Q}(\sqrt{p})$, then Chapman \cite{Chapman} proved that 
    $$\det C_1(t)=\begin{cases}
    	(-1)^{\frac{p-1}{4}}2^{\frac{p-1}{2}}(-a_pt+b_p) & \mbox{if}\ p\equiv 1\pmod4,\\
    	-2^{\frac{p-1}{2}}t                              & \mbox{if}\ p\equiv 3\pmod4,
    \end{cases}$$
    where $a_p,b_p\in\mathbb{Q}$ are defined by the equality 
    $$\varepsilon_p^{h_p}=a_p+b_p\sqrt{p}.$$
    
    In 2019, Sun \cite{Sun} initiated the study of determiants involving the Legendre symbol and binary quadratic forms. For example, Sun considered the determinant
    $$\det S_p:=\det\left[\left(\frac{i^2+j^2}{p}\right)\right]_{1\le i,j\le (p-1)/2}.$$
    Sun \cite[Theorem 1.2]{Sun} showed that $-\det S_p$ is always a quadratic residue modulo $p$. Readers may refer to \cite{krachun, W21} for the recent works on this topic. 
    
    Recently, Sun \cite{Sun-Conj} posed many interesting conjectures on determinants related to the Legendre symbol. For example, Sun \cite[Conjecture 1.1]{Sun-Conj} posed the following conjecture.
    
    \begin{conjecture}[Sun]
    	Let $p\equiv3\pmod 4$ be a prime. Then 
    	$$\det\left[t+\left(\frac{i^2+j^2}{p}\right)+\left(\frac{i^2-j^2}{p}\right)\right]_{1\le i,j\le (p-1)/2}=\left(\frac{p-1}{2}t-1\right)p^{\frac{p-3}{4}}.$$
    \end{conjecture}
    
    Motivated by the above results, in this paper, we will study some determinants involving the quadratic multiplicative character of a finite field. We first introduce some notations. 
    
    Let $q=p^s$ be an odd prime power with $p$ prime and $s\in\mathbb{Z}^+$ and let $\mathbb{F}_q$ be the finite field of $q$ elements. Let $\mathbb{F}_q^{\times}$ be the cyclic group of all nonzero elements of $\mathbb{F}_q$. For any positive integer $k\mid q-1$, let 
    $$D_k:=\{a_1,a_2,\cdots,a_{(q-1)/k}\}=\{x^k:\ x\in\mathbb{F}_q^{\times}\}$$
    be the subgroup of all nonzero $k$-th powers in $\mathbb{F}_q$.
    
    Let $\widehat{\mathbb{F}_q^{\times}}$ be the cyclic group of all multiplicative characters of $\mathbb{F}_q$. Throughout this paper, for any $\psi\in\widehat{\mathbb{F}_q^{\times}}$, we extend $\psi$ to $\mathbb{F}_q$ by setting $\psi(0)=0$. Also, if $2\nmid q$, then we use the symbol $\phi$ to denote the unique quadratic multiplicative character of $\mathbb{F}_q$, i.e., 
    \begin{equation*}
    	\phi(x)=\begin{cases}
    		1  & \mbox{if}\ x\in D_2,\\
    		0  & \mbox{if}\ x=0,\\
    		-1 & \mbox{otherwise.}
    	\end{cases}
    \end{equation*}
    
    Inspired by the above results, in this paper, we define the matrix $A_k(t)$ by 
    $$A_k(t):=\left[t+\phi(a_i+a_j)+\phi(a_i-a_j)\right]_{1\le i,j\le (q-1)/k}.$$
    Also, the integers $c_k$ and $d_k$, which are related to number of $\mathbb{F}_q$-rational points of hyperelliptic curves over $\mathbb{F}_q$, are defined by 
    \begin{equation}\label{Eq. definition of ak}
    	\left|\{\infty\}\cup\{(x,y)\in\mathbb{F}_q\times\mathbb{F}_q:\ y^2=x^k+1\}\right|=q+1-c_k
    \end{equation}
    and 
    \begin{equation}\label{Eq. definition of bk}
    	\left|\{\infty\}\cup\{(x,y)\in\mathbb{F}_q\times\mathbb{F}_q:\ y^2=x^k-1\}\right|=q+1-d_k.
    \end{equation}
    
    Now we state our main results of this paper.
    
    \begin{theorem}\label{Thm. Aq(k)}
    	 Let $q=p^s$ be an odd prime power with $p$ prime and $s\in\mathbb{Z}^+$. Then for any positive integer $k\mid q-1$, the following results hold.
    	 
    	 {\rm (i)} Suppose $q\equiv 1\pmod {2k}$. Then $\det A_k(t)=0$. In particular, in the case $q\equiv1\pmod4$, we have $\det A_2(t)=0$. 
    	 
    	 {\rm (ii)} If $q\equiv 3\pmod4$, then 
    	 $$\det A_2(t)=\left(\frac{q-1}{2}t-1\right)q^{\frac{q-3}{4}}.$$
    	 
    	 {\rm (iii)} Suppose $q\equiv 1\pmod4$ and $q\not\equiv1\pmod{2k}$. Then there is an integer $u_k$ such that 
    	 $$\det A_k(t)=\left(\frac{q-1}{k}t-\frac{1}{k}(c_k+d_k+2)\right)\cdot u_k^2.$$
    \end{theorem}
    
    \begin{remark}
    	(i) Theorem \ref{Thm. Aq(k)}(i) generalizes the result \cite[Theorem 1.1]{Sun-Conj} to an arbitrary finite field with odd characteristic. In the case $q=p$ is an odd prime, Theorem \ref{Thm. Aq(k)}(ii) confirms the above conjecture \cite[Conjecture 1.1]{Sun-Conj} posed by Zhi-Wei Sun.
    	
    	(ii) For any $3\le k<q-1$ with $k\mid q-1$ and $q-1\not\equiv 0\pmod{2k}$, we can also obtain the explicit value of $\det A_k(t)$. However, finding a simple expression of $\det A_k(t)$ seems very difficult.
    \end{remark}
    
     We will prove our main results in Section 2.

    \section{Proof of Theorem \ref{Thm. Aq(k)}}

    Throughout this section, we let $\chi$ be a generator of $\widehat{\mathbb{F}_q^{\times}}$. Also, for any $\chi^i,\chi^j\in\widehat{\mathbb{F}_q^{\times}}$, the Jacobi sum of  $\chi^i$ and $\chi^j$ is defined by 
    $$J(\chi^i,\chi^j)=\sum_{x\in\mathbb{F}_q}\chi^i(x)\chi^j(1-x).$$
    
   We begin with a known result in linear algebra.
   
   \begin{lemma}\label{Lem. eigenvalues}
   		Let $n$ be a positive integer and let $M$ be an $n\times n$ complex matrix. Let $\lambda_1,\cdots,\lambda_n\in\mathbb{C}$, and let $\v_1,\cdots,\v_n\in\mathbb{C}^n$ be column vectors. Suppose that 
   	$$M\v_i=\lambda_i\v_i$$
   	for each $1\le i\le n$ and that the vectors $\v_1,\cdots,\v_n$ are linearly independent over $\mathbb{C}$. Then $\lambda_1,\cdots,\lambda_n$ are exactly all the eigenvalues of $M$ (counting multiplicity). 
   \end{lemma}
	
	Before the proof of our main results, we first introduce the definition of circulant matrices. Let $R$ be a commutative ring and let $b_0,b_1,\cdots,b_{n-1}\in R$. Then the circulant matrix of the tuple $(b_0,b_1,\cdots,b_{n-1})$ is defined by 
	$$C(b_0,b_1,\cdots,b_{n-1}):=[b_{i-j}]_{0\le i,j\le n-1},$$
	where the indices are cyclic modulo $n$. 
	
	The second author \cite[Lemma 3.4]{W21} proved the following result.
	
	\begin{lemma}\label{Lem. Wuffa}
		Let $n\ge1$ be an odd integer. Let $R$ be a commutative ring and let $b_0,\cdots,b_{n-1}\in R$ such that 
		$$b_i=b_{n-i}$$
		for any $1\le i\le n-1$. Then there is an element $u\in R$ such that 
		$$\det C(b_0,b_1,\cdots,b_{n-1})=\left(\sum_{i=0}^{n-1}b_i\right)u^2.$$
	\end{lemma}
	
	Now we are in a position to prove our first result. For simplicity, we set $n=(q-1)/k$.
	
	{\noindent{\bf Proof of Theorem \ref{Thm. Aq(k)}.}}  (i) Suppose $q-1\equiv 0\pmod{2k}$. Let $\xi_{2k}\in\mathbb{F}_q$ be a primitive $2k$-th root of unity. Then $-1=\xi_{2k}^k\in D_k$. Thus, for any $1\le j\le n$ there exists an integer $1\le j'\le n$ such that $a_{j'}=-a_j$ and $j\neq j'$. This implies that the $j$-th column of $A_k(t)$ is the same as the $j'$-th column of $A_k(t)$ and hence $\det A_k(t)=0$. 
	
	(ii) Suppose now $q-1\not\equiv 0\pmod {2k}$. Then clearly $k$ is even.	For any integers $0\le m\le n-1$ and $1\le i\le n$, we have 
	\begin{align*}
		 &\sum_{1\le j\le n}\left(\phi(a_i+a_j)+\phi(a_i-a_j)\right)\chi^m(a_j)\\
		=&\sum_{1\le j\le n}\left(\phi\left(1+\frac{a_j}{a_i}\right)+\phi\left(1-\frac{a_j}{a_i}\right)\right)\chi^m\left(\frac{a_j}{a_i}\right)\chi^m(a_i)\\
		=&\sum_{1\le j\le n}\left(\phi(1+a_j)+\phi(1-a_j)\right)\chi^m(a_j)\chi^m(a_i).
	\end{align*}
	Let 
	$$\v_m=\left(\chi^m(a_1),\chi^m(a_2),\cdots,\chi^m(a_n)\right)^T,$$
	and let 
	$$\lambda_m=\sum_{1\le j\le n}\left(\phi(1+a_j)+\phi(1-a_j)\right)\chi^m(a_j).$$
	Then by the above results, for any $0\le m\le n-1$ we obtain 
	\begin{equation*}
		A_k(0)\v_m=\lambda_m\v_m.
	\end{equation*}
	Since 
	$$\det \left[\chi^i(a_j)\right]_{0\le i\le n-1,1\le j\le n}=\prod_{1\le i<j\le n}\left(\chi(a_j)-\chi(a_i)\right)\neq 0,$$
	the vectors $\v_0,\cdots,\v_{n-1}$ are linearly independent over $\mathbb{C}$, and hence by Lemma \ref{Lem. eigenvalues} the numbers $\lambda_0,\cdots,\lambda_{n-1}$ are exactly all the eigenvalues of $A_k(0)$. 
	
	Now let $k=2$. Then clearly $q\equiv 3\pmod 4$ and $n$ is odd in this case. We first evaluate $\det A_2(0)$. By the above, we have 
	\begin{equation}\label{Eq. a1 in the proof of Thm. Aq(k)}
		\det A_2(0)=\lambda_0\prod_{1\le m\le n-1}\lambda_m=\lambda_0\prod_{1\le m\le (n-1)/2}\left|\lambda_{2m}\right|^2.
	\end{equation}
	The last equality follows from $\overline{\lambda_m}=\lambda_{n-m}$ for $1\le m\le n-1$. For $\lambda_0$ we have 
	\begin{align}\label{Eq. a2 in the proof of Thm. Aq(k)}
		  \lambda_0
		&=\sum_{1\le j\le n}\left(\phi(1+a_j)+\phi(1-a_j)\right) \notag \\
		&=\frac{1}{2}\sum_{x\in\mathbb{F}_q^{\times}}\phi(1+x^2)-\frac{1}{2}\sum_{x\in\mathbb{F}_q^{\times}}\phi(x^2-1) \notag\\
		&=-1.
	\end{align}
    The last equality follows from 
    $$\sum_{x\in\mathbb{F}_q}\phi(x^2\pm 1)=-1.$$
    For $\lambda_{2m}$ with $1\le m\le (n-1)/2$, one can verify that 
    \begin{align}\label{Eq. a3 in the proof of Thm. Aq(k)}
    	  \lambda_{2m}
    	&=\sum_{1\le j\le n}\left(\phi(1+a_j)+\phi(1-a_j)\right)\chi^{2m}(a_j) \notag \\
    	&=\frac{1}{2}\sum_{x\in\mathbb{F}_q}\phi(1+x^2)\chi^{2m}(x^2)+\frac{1}{2}\sum_{x\in\mathbb{F}_q}\phi(1-x^2)\chi^{2m}(-x^2) \notag \\
    	&=\sum_{x\in\mathbb{F}_q}\phi(1+x)\chi^{2m}(x) \notag \\
    	&=\sum_{x\in\mathbb{F}_q}\phi(1+x)\chi^{2m}(-x) \notag \\
    	&=J(\phi,\chi^{2m}).
    \end{align}
    
	Combining (\ref{Eq. a2 in the proof of Thm. Aq(k)}) and (\ref{Eq. a3 in the proof of Thm. Aq(k)}) with (\ref{Eq. a1 in the proof of Thm. Aq(k)}), we obtain 
	$$\det A_2(0)=-\prod_{1\le m\le (n-1)/2}\left|J(\phi,\chi^{2m})\right|^2=-q^{\frac{q-3}{4}}.$$
	
	Now we turn to $\det A_2(t)$. By (\ref{Eq. a2 in the proof of Thm. Aq(k)}) for any $1\le j\le n$ we have 
	\begin{align*}
			 &\sum_{1\le i\le n}\left(t+\phi(a_i+a_j)+\phi(a_i-a_j)\right)\\
			=& nt+\sum_{1\le i\le n}\left(\phi(1+a_j/a_i)+\phi(1-a_j/a_i)\right)\\
			=& nt++\sum_{1\le i\le n}\left(\phi(1+a_i)+\phi(1-a_i)\right)\\
			=& nt-1.
	\end{align*}
    This implies that $(nt-1)\mid \det A_2(t)$. Noting that $\det A_2(t)\in \mathbb{Z}[t]$ with degree $\le 1$, we obtain 
    $$\det A_2(t)=-\det A_2(0)\cdot(nt-1)=q^{\frac{q-3}{4}}\left(\frac{q-1}{2}t-1\right).$$
	
	(iii) Suppose $q\equiv 1\pmod 4$ and $q\not\equiv 1\pmod {2k}$. Clearly $k\equiv 0\pmod 2$ in this case. Let $g\in\mathbb{F}_q$ be a generator of the cyclic group $\mathbb{F}_q^{\times}$. Then one can verify that 
	\begin{align*}
		 \det A_k(t)
		 &=\det \left[t+\phi(a_i+a_j)+\phi(a_i-a_j)\right]_{1\le i,j\le n}\\
	     &=\det \left[t+\phi(g^{k(i-j)}+1)+\phi(g^{k(i-j)}-1)\right]_{0\le i,j\le n-1}.
	\end{align*}
	For $0\le i\le n-1$, let 
	$$b_i=t+\phi(g^{ki}+1)+\phi(g^{ki}-1).$$ 
	Then one can easily verify that 
	$$\det A_k(t)=\det C(b_0,b_1,\cdots,b_{n-1}),$$
	and that $b_i=b_{n-i}$ for any $1\le i\le n-1$. Now applying Lemma \ref{Lem. Wuffa} we see that there is an element $u_k\in\mathbb{Z}[t]$ such that 
	$$\det A_k(t)=\left(\sum_{i=0}^{n-1}b_i\right)\cdot u_k^2.$$
	One can verify that 
	\begin{align*}
		 \sum_{i=0}^{n-1}b_i
	   &=nt+\sum_{1\le j\le n}\left(\phi(a_i+1)+\phi(a_i-1)\right)\\
	   &=nt+\frac{1}{k}\sum_{x\in\mathbb{F}_q^{\times}}\left(\phi(x^k+1)+\phi(x^k-1)\right)\\
	   &=nt-\frac{1}{k}(c_k+d_k+2),
	\end{align*}
    Where $c_k$ and $d_k$ are defined by (\ref{Eq. definition of ak}) and (\ref{Eq. definition of bk}), and the last equality follows from 
    $$\sum_{x\in\mathbb{F}_q^{\times}}\phi(x^k+1)=-c_k-1$$
    and 
    $$\sum_{x\in\mathbb{F}_q^{\times}}\phi(x^k-1)=-d_k-1.$$
    
	As $\det A_k(t)\in\mathbb{Z}[t]$ with degree $\le 1$, by the above we see that $u_k\in\mathbb{Z}$. Hence 
	$$\det A_k(t)=\left(\frac{q-1}{k}t-\frac{1}{k}(c_k+d_k+2)\right)\cdot u_k^2.$$
	
	In view of the above, we have completed the proof of Theorem \ref{Thm. Aq(k)}.\qed

\end{document}